\documentstyle[11pt,a4wide]{article}
\def\?{\bigskip$$\vdots$$\bigskip}
\begin{document}
\bibliographystyle{plain}

Syst\`emes Dynamiques / {\sl Dynamical Systems}.

\begin{center}
{\bf 
Perturbation d'un hamiltonien partiellement hyperbolique
}
\end{center}
\begin{center}
Patrick {\sc Bernard}
\end{center}
\begin{small}
{\bf R\'esum\'e:}
\begin{it}
La diffusion d'Arnold pour les syst\`emes hamiltoniens presque int\'egrables
$H(q,p)=H_0(p)+\mu h(q,p)$ se fait en temps exponentiellement grand
en $\mu$. On \'etudie ici un syst\`eme initialement hyperbolique qui
~admet une diffusion en temps polynomial.
\end{it}
\begin{center}
\bf
Perturbation of a partially hyperbolic hamiltonian system
\end{center} 
{\bf Abstract:}
\begin{it}
Arnold's diffusion in quasi integrable hamiltonian systems
$H(q,p)=H_0(p)+\mu h(q,p)$ occurs in exponentially large time.
We study an initially hyperbolic
system which admits diffusion in polynomial time.
\end{it}
\end{small}

{\bf Abridged english version}

 In spite of KAM theorem, quasi integrable hamiltonian systems 
%($ H(p,q) = H_0(p)+\epsilon H_1(p,q),$
%  $ (p,q)\in {\mbox{{\bf R}}}^n \times {\mbox{{\bf T}}}^n $) 
need not be stable.
In \cite{Arnold}, V.I.  Arnold has pointed out an unstable behavior for
systems with more than three degrees of freedom.
Such systems can admit orbits satisfying $|p(t_1)-p(t_0)|\ge 1$.
Arnold's example can be expressed as  a lagrangian depending on two parameters:
$$ { \cal L} _{\epsilon,\mu}(q,\dot q , Q,\dot Q , t) =
\frac 1 2 \dot Q^2 +\frac 1 2 \dot q^2+
\epsilon (1-\cos q)
+\epsilon\mu (1-\cos q)(\cos Q +\cos t). $$
The first perturbation ${\cal L}_{\epsilon ,0}$ preserves a family of 
invariant tori {\bf T}$_{\omega}$ ($\dot q = q = \dot Q-\omega =0$)
and makes them partially hyperbolic. The second perturbation creates 
connections between these '' whiskered`` tori. There now exists unstable
orbits sliding along these connections. In a recent paper
\cite{Bessi} U. Bessi gives a variationnal proof of this result. He 
obtains an estimate of the smallest time $t_d $ for which 
$|\dot Q(t_d)-\dot Q(0)|\ge 1$ : for $\epsilon $ small enough and  
$\mu \le e^{-a/\epsilon}$, $t_d \le e^{c/\sqrt{\epsilon}}$. 
We recall that this time must be exponentially high in  $\epsilon $ as was 
proved in  \cite{Lochak}. In both these methods partial hyperbolicity 
of {\bf T}$_{\omega}$ is essential. That's why it seems useful to study
initially hyperbolic systems for which diffusion, though  much quicker and 
easier to point out, is not completely understood.
Following a suggestion by P. Lochak, I will study the example 
$$ {\cal L} _{\mu}(q,\dot q , Q,\dot Q , t) =
  \frac 1 2 \dot Q^2 +\frac 1 2 \dot q^2+
  (1-\cos q)
  +\mu (1-\cos q)f(Q,q,t)
$$
where f is $2\pi$ periodic in each variable and, using Bessi's method, I will
prove the folowing result:

{\sc Theorem :}
{\it
For $\mu$ small enough, and f satisfying condition 1, $\cal L_{\mu }$ admits 
orbits satisfying:
 $$ \dot Q(0) \le 1\ , \dot Q (t_d) \ge 2
 \mbox{~~with~~} t_d \le C \mu^{-2}  \ . $$
}

Condition 1 will be expressed later. The polynomial time obtained there is 
much lower than the time obtained for initially integrable systems.

We first investigate the nonperturbed system. Every torus 
$ \mbox{\bf T} _\omega : q=\dot q= \dot Q-\omega =0$ is invariant and admits 
a two parameters family of homoclinic orbits :
$ q(t)=q_0(t-T_0)\ ,\ Q(t)=\omega(t-T_0)+Q_0 $
where $ q_0 (t)=4\arctan e ^t $.
Let 
$$ {\cal A}_{\omega} (Q_0 ,T_0) = -\int  _{-\infty } ^{+\infty  }
(1-\cos q_0 (t-T_0))\,f(Q_0+\omega(t-T_0),q_0(t-T_0),t) \,dt $$
be the perturbation part of the action along one of these homoclinic orbits.
The function $ {\cal A}_{\omega} $ is periodic, so it must have a minimum. 
Condition 1 says that this minimum is nondegenerate.

{\sc Condition 1:}
There exists a nonnegative real number A which satisfies : For every  
$\omega \in \left[ \frac 1 2 ,\frac 5 2 \right]$, we can find 
 $(Q_\omega,T_\omega) \in \mbox{\bf R}^2$ 
and a box ${\cal B}(\omega)=[T_\omega -s_\omega,T_\omega+s_\omega]$
with $|r_\omega|<\pi$ and 
 $|s_\omega|<\pi$
such that ${\cal A}(\omega,Q,T)
   \ge {\cal A}(\omega,Q_\omega,T_\omega)+A$ for any
$ (Q,T)\in \partial {\cal B}(\omega).$ 

There are two steps in the proof. We first notice that there exists
trajectories of the perturbed system that start at the top of an unperturbed 
homoclinic orbit, follow it down to an invariant torus, stay a long 
time around that torus, and then leave it along another unperturbed 
homoclinic orbit.
We then glue orbits of that kind together to obtain connections between 
close to each other invariant tori and we can move step by step
along the chain of invariant tori following  these connections. 
The 
time of diffusion is then obtained as the quotient between the time 
$T=1/\mu$
of one 
of these transitions and the distance $\mu$ between two tori 
that can be joined
by such a transition.

\begin{center} 
------
\end{center}

Pour un syst\`eme hamiltonien presque 
int\'egrable 
($ H(p,q) = H_0(p)+\epsilon H_1(p,q),
   (p,q)\in {\mbox{{\bf R}}}^n \times {\mbox{{\bf T}}}^n $),
le th\'eor\`eme KAM donne l'existence
d'un grand nombre de tores invariants et donc d'un grand nombre de 
trajectoires stables, induisant une petite variation de l'action $p$.
Ceci ne suffit pourtant pas \`a garantir la stabilit\'e du syst\`eme.
Dans \cite{Arnold},V.I.  Arnold a mis en \'evidence  un m\'ecanisme 
d'instabilit\'e pour les syst\`emes presque int\'egrables ayant au moins trois 
degr\'es de libert\'e. Il peut exister, pour de tels syst\`emes, des
trajectoires instables, induisant une grande d\'erive de l'action 
c'est \`a dire telles que
$|p(t_1)-p(t_0)|\ge 1$. L'exemple d'Arnold est une 
perturbation \`a deux param\`etres d'un lagrangien int\'egrable :
$$ { \cal L} _{\epsilon,\mu}(q,\dot q , Q,\dot Q , t) =
\frac 1 2 \dot Q^2 +\frac 1 2 \dot q^2+
\epsilon (1-\cos q)
+\epsilon\mu (1-\cos q)(\cos Q +\cos t). $$
La premi\`ere perturbation ${\cal L}_{\epsilon ,0}$ conserve une famille de 
tores invariants {\bf T}$_{\omega}$ ($\dot q = q = \dot Q-\omega =0$),
qu'elle rend partiellement hyperboliques.
La deuxi\`eme perturbation 
cr\'ee des connections entre tores voisins dans cette famille, qui reste
invariante. 
Il existe alors des trajectoires qui
longent ces connections et induisent 
une d\'erive du terme $\dot Q$ de l'action.
Dans un article r\'ecent \cite{Bessi}, U. Bessi a \'etudi\'e 
cet exemple par une
m\'ethode variationnelle. Il obtient une majoration du temps de diffusion 
$t_d $ tel que $|\dot Q(t_d)-\dot Q(0)|\ge 1$: pour $\epsilon $
petit et $\mu$ exponentiellement petit devant $\epsilon$ 
($\mu \le e^{-a/\epsilon}$), il existe des trajectoires diffusantes  qui
v\'erifient $t_d \le e^{c/\sqrt{\epsilon}}$. 
Remarquons que ce temps a \'et\'e 
minor\'e dans \cite{Lochak} par 
$t_{d}\ge \exp(\epsilon^{-1/4})$, la diffusion dans ces syst\`emes est 
exponentiellement lente. 
De nombreux probl\`emes se posent lorsqu'on essaie de g\'en\'eraliser ces 
m\'ethodes (voir \cite{Bar} pour un r\'esum\'e des id\'ees actuelles sur 
ce sujet), mais il appara\^\i t que l'introduction d'hyperbolicit\'e partielle 
par une premi\`ere perturbation est une \'etape importante. 
Il semble donc utile
d'\'etudier des syst\`emes directement hyperboliques, pour lesquels 
la diffusion, quoique  beaucoup plus claire et beaucoup plus 
rapide \cite{Bar}, n'est pas encore comprise.  
\`A la suite d'une suggestion de P. Lochak, je me suis donc interess\'e \`a
l'exemple
$$ {\cal L} _{\mu}(q,\dot q , Q,\dot Q , t) =
  \frac 1 2 \dot Q^2 +\frac 1 2 \dot q^2+
  (1-\cos q)
  +\mu (1-\cos q)f(Q,q,t), 
$$
o\`u $f$ est une fonction $2\pi$ p\'eriodique en chaque variable.
Par la m\'ethode de Bessi, on montre le r\'esultat suivant:

{\sc Th\'eor\`eme :}
{\it
Pour $\mu$ suffisamment petit,
le syst\`eme de lagrangien $\cal L_{\mu }$
, o\`u $f$ v\'erifie la condition 1 ci-dessous
, admet des trajectoires qui v\'erifient:
$$ \dot Q(0) \le 1\ , \dot Q (t_d) \ge 2
 \mbox{~~avec~~} t_d \le C \mu^{-2}  \ . $$
}

La condition 1, qui sera explicit\'ee dans la suite, ne repr\'esente qu'une 
contrainte minime, elle est en particulier v\'erifi\'ee par le cas $\epsilon=1$
de l'exemple d'Arnold. 
Le temps de diffusion obtenu, polynomial,  est effectivement beaucoup plus 
faible que celui d'un syst\`eme initialement int\'egrable.  

{\sc 1- Syst\`eme non perturb\'e:}
Dans le cas non perturb\'e on a affaire au 
produit non coupl\'e d'un oscillateur libre et
d'un pendule pesant.
Ce syst\`eme admet une famille \`a un param\`etre de tores
invariants partiellement hyperboliques
$ \mbox{\bf T} _\omega : q=\dot q= \dot Q-\omega =0.$
L'\'equation de la s\'eparatrice du pendule valant $\pi$ au temps $T_0$ est
$ q(t)=q_0 (t-T_0 ) \ ,$
o\`u
$ q_0 (t)=4\arctan e ^t $,
$\dot q_0 (t)=2 (\cosh t)^{-1}  $.
Chaque tore {\bf T}$_\omega$ admet une famille \`a deux param\`etres
de trajectoires homoclines:
$ q(t)=q_0(t-T_0)\ ,\ Q(t)=\omega(t-T_0)+Q_0.$
On appellera par la suite s\'eparatrices ces trajectoires homoclines non
perturb\'ees.
On consid\`erera la fonctionnelle suivante, qui est 
le terme de perturbation de l'action le long d'une s\'eparatrice:
$$ {\cal A}_{\omega} (Q_0 ,T_0) = -\int  _{-\infty } ^{+\infty  }
(1-\cos q_0 (t-T_0))\,f(Q_0+\omega(t-T_0),q_0(t-T_0),t) \,dt. $$
La fonction $\cal A _\omega$ est p\'eriodique, et admet un minimum.
On supposera  dans la suite v\'ertifi\'ee la condition 1 suivante, 
qui est une condition de non 
d\'eg\'enerescence de ce minimum:

{\sc Condition 1:}
Il existe un r\'eel positif A tel que
pour tout $\omega \in \left[ \frac 1 2 ,\frac 5 2 \right]$,
il existe un point 
 $(Q_\omega,T_\omega) \in \mbox{\bf R}^2$ 
et une boite  
 ${\cal B}(\omega)=[T_\omega -s_\omega,T_\omega+s_\omega]
\times [Q_\omega-r_\omega,Q_\omega+r_\omega]$ avec $|r_\omega|<\pi$ et
 $|s_\omega|<\pi$,
tels que l'on ait 
${\cal A}(\omega,Q,T)
   \ge {\cal A}(\omega,Q_\omega,T_\omega)+A$ pour tout  
$ (Q,T)\in \partial {\cal B}(\omega).$

2- Le lemme 1 ci-dessous donne l'existence d'une trajectoire qui longe une
s\'eparatrice depuis
son apog\'ee $(q=\pi)$ jusqu'\`a un tore invariant,
puis reste le long du tore,
puis repart le long d'une
autre s\'eparatrice. On admettra ce lemme, dont 
la d\'emonstration demande quelques calculs,
que l'on peut trouver dans \cite{Bessi}.
On pose  $T=A_0/\mu$, avec $A_0>1$ et $l=|\log \mu|/6$. 

{\sc Lemme 1}
{\it
Soit $\mu$ petit et
$(Q_0,Q_1,T_0,T_1)\in \mbox{\bf R}^4$.
Il existe une unique trajectoire minimisante
$(Q,q,\dot Q,\dot q)$ du syst\`eme perturb\'e v\'erifiant les conditions
aux extr\'emit\'es
$q(T_0)=-\pi ,q(T_1)=\pi ,Q(T_0)=Q_0,Q(T_1)=Q_1$.
De plus, il existe  $\alpha>0$ tel que si $T\le T_1-T_0 $
la trajectoire ainsi d\'efinie soit $\mu^\alpha$-proche en topologie $C^0$
de la
s\'eparatrice $(\omega,T_0,Q_0)$ sur $[T_0,T_0+l]$, puis du tore
{\bf T}$_\omega$ sur $[T_0+l,T_1-l]$ puis de la s\'eparatrice
$(\omega,T_1,Q_1)$, o\`u
$
\omega  =(Q_1-Q_0)/(T_1-T_0).
$
}

On associe ainsi, \`a un tore {\bf T}$_\omega$ et \`a deux s\'eparatrices de 
ce tore, une solution particuli\`ere, une {\it boucle}, qui v\'erifie les 
propri\'et\'es d\'ecrites ci-dessus. On a une famille
 \`a trois param\`etres de boucles 
venant tourner autour d'un tore fix\'e ,
puisque choisir un tore revient \`a imposer la relation 
$
\omega  =(Q_1-Q_0)/(T_1-T_0)$ entre les conditions aux bords.
En recollant deux boucles, on obtient une 
{\it transition simple}, c'est \`a dire 
une trajectoire qui joint un voisinage d'un 
premier tore {\bf T}$_{\omega_1}$ \`a un voisinage d'un second tore
{\bf T}$_{\omega_2}$. C'est l'objet de la proposition suivante, qui sera 
d\'emontr\'ee par la suite.   
  
{\sc Proposition}:
\begin{it}  
Pour $\mu$ suffisamment petit, pour tous $\omega_1$
et $\omega_2$ tels que $1/2<\omega_1<\omega_2<5/2$ et 
$0 \le \omega_2 - \omega_1\le \mu$,
pour tous $(T_{-1},Q_{-1},T_1,Q_1)$ tels que 
$|(Q_1-Q_{-1})/(T_1-T_{-1})-\omega_1|\le \mu$ et $T_1-T_{-1}\ge  2T $,  
il existe $(T_0,Q_0)$, tel que 
$|(Q_1-Q_0)/(T_1-T_0)-\omega_2|\le \mu$ et 
$|(Q_0-Q_{-1})/(T_0-T_{-1})-\omega_1|\le \mu$, et
tel que la juxtaposition des boucles
d'extr\'emit\'es $(T_{-1},Q_{-1},T_0,Q_0)$ et $(T_0,Q_0,T_1,Q_1)$
forme, sur $[T_{-1},T_1]$, une solution du syst\`eme.
\end{it}

Le th\'eor\`eme en d\'ecoule simplement, en effet
on obtient ainsi non seulement une transition simple joignant 
{\bf T}$_{\omega_1}$ et {\bf T}$_{\omega_2}$, mais une famille 
de telles transitions, les param\`etres \'etant 
ici les valeurs aux extr\'emit\'es  $(T_{-1},Q_{-1},T_1,Q_1)$. 
On peut donc reproduire le processus de recollement pour obtenir des
transitions doubles (longeant trois tores) puis triples, etc \dots
On peut ainsi construire des transitions multiples joignant pas \`a 
pas deux tores arbitrairement \'eloign\'es 
{\bf T}$_{\omega_1}$ et {\bf T}$_{\omega_2}$, si $1/2< \omega_1 \le
\omega_2<5/2$.
Le temps de diffusion 
appara\^\i t comme le quotient du temps $T$ d'une transition simple 
par l'\'ecart entre deux tores que l'on peut joindre par une telle 
transition. Ceci d\'emontre le th\'eor\`eme puisque  $T=A_0/\mu$.

{\sc 3- D\'emonstration de la proposition}:
Soit $H$ l'espace des courbes de $H^1([T_{-1},T_1],\mbox{\bf R}^2)$
v\'erifiant les conditions
aux extr\'emit\'es suivantes: $q(T_{-1})=-\pi$, $q(T_1)=3\pi$, $Q(T_{-1})=
Q_{-1}$, $Q(T_1)=Q_1$. On va construire, par une m\'ethode variationnelle,
une transition simple dans cet espace. Soit $\iota : \mbox{\bf R}^2
\longrightarrow H$ l'application qui, \`a $(T_0,Q_0)\in \mbox{\bf R}^2$, 
associe la courbe obtenue par recollement des deux boucles  
$(T_{-1},Q_{-1},T_0,Q_0)$ et $(T_0,Q_0,T_1,Q_1)$.
On d\'efinit sur $H$ la fonctionnelle d'action
$\cal F$
et on notera ${\cal F}_Y$ aussi bien sa restriction \`a
{\bf Y} $=\iota(\mbox{\bf R}^2)$ que ${\cal F}\circ \iota$. 

{\sc Lemme 2}:
\begin{it}
Un minimum local de ${\cal F}_Y$ est  une solution du syst\`eme 
sur $[T_{-1},T_1]$.
\end{it}

Pour d\'emontrer ce lemme,
il suffit de montrer que tout minimum local sur {\bf Y} est un minimum local
sur $H$. Pour ceci consid\'erons l'application 
$G:H \longrightarrow \mbox{\bf R}^2$, $(Q(t),q(t)) \longmapsto (a,b)$, o\`u
$a=\inf \{t:q(t)>\pi\}$ et $b=q(a)$. 
On notera aussi $G$ l'application $\iota \circ G$.
Les boucles \'etant des trajectoires minimisantes,
${\cal F}(\gamma)\ge {\cal F}(G(\gamma))$. 
Par ailleurs les boucles, proches des s\'eparatrices, v\'erifient 
$\dot q(T_0)>0$ et $q(t)\neq q(T_0)$ pour $t\neq T_0$ . 
Donc $G$ est continue au voisinage de {\bf Y}.
Soit alors $\gamma_0$ un minimum local de ${\cal F}_Y$, et soit
$\gamma \in H$. Si $\gamma$ est proche de $\gamma_0$, alors 
$G(\gamma)$ est proche de $\gamma_0=G(\gamma_0)$, donc
${\cal F}(\gamma)\ge {\cal F}(G(\gamma)) \ge {\cal F}(\gamma_0)$;
$\gamma_0$ est donc un minimum local de $\cal F$ sur $H$, ceci d\'emontre
le lemme 2.  

On est donc ramen\'e \`a chercher un minimum local bien plac\'e de 
${\cal F}_Y$ pour d\'emontrer la proposition.
Les boucles longeant les s\'eparatrices auxquelles elles sont associ\'ees,
on peut approcher leurs actions par celles de ces s\'eparatrices. En fait,
l'action d'une s\'eparatrice n'est pas d\'efinie. Par contre le terme
de perturbation $\cal P$ de cette action est bien d\'efini. On va  
chercher dans un premier temps \`a minimiser $\cal P$, puis on montrera
que ses variations sont pr\'edominantes dans la zone consid\'er\'ee, et donc 
que $\cal F$ admet un minimum local au voisinage du minimum local de $\cal P$.
Etudions $\cal P$ le long d'une boucle:
avant d'atteindre le tore {\bf T}$_{\omega_1}$, il y a une partie
de courbe qui longe une s\'eparatrice fix\'ee (ne d\'ependant pas de 
$(T_0,Q_0)$), puis une partie de courbe entre {\bf T}$_{\omega_1}$ et 
{\bf T}$_{\omega_2}$ longeant la s\'eparatrice $(\omega_1,T_0,Q_0)$,
puis une autre partie longeant un morceau fix\'e de s\'eparatrice.
Par cons\'equent ${\cal P}(T_0,Q_0)\approx C + {\cal A}(\omega_1,T_0,Q_0)$, 
et il est naturel de 
chercher un minimum local de $\cal P$ au voisinage du minimum 
$(T_{\omega_1},Q_{\omega_1})$ de ${\cal A}_{\omega_1}$.

Plus pr\'ecis\'ement,
soit $B\subset${\bf R}$^2$ tel que $B = {\cal B}(T_{\omega_1},Q_{\omega_1})$
mod$[2\pi]$ et tel que pour tout  $ (T_0,Q_0)\in B$ on ait $ 
|(Q_1-Q_0)/(T_1-T_0)-\omega_2|\le \mu$
et  $|(Q_0-Q_{-1})/(T_0-T_{-1})-\omega_1|\le \mu$. 
$\cal F _Y$  admet un minimum sur $B$,
qui est compact. Si ce minimum est atteint \`a l'int\'erieur de $B$, il 
s'agit d'un minimum local de ${\cal F}_Y$ et la proposition d\'ecoule du 
lemme 2. Montrons que
${\cal F}(\partial B)\ge {\cal F}(T_{\omega_1},Q_{\omega_1})$;
pour ceci, posons ${\cal F=F}_0+\mu {\cal P}$. 
On peut appliquer le lemme 1 au cas 
$\mu = 0$ et  on obtient des boucles non perturb\'ees. On notera
$\gamma_{\omega_1}$ la trajectoire obtenue par recollement des boucles non
perturb\'ees entre $(T_{-1},Q_{-1},T_{\omega_1},Q_{\omega_1})$ et 
$(T_{\omega_1},Q_{\omega_1},T_1,Q_1)$.
Remarquons que ${\cal F}(\gamma_{\omega_1})\ge 
{\cal F}(T_{\omega_1},Q_{\omega_1})$.
Il est utile de consid\'erer cette trajectoire car il est plus facile d'estimer
le terme de perturbation de l'action le long d'une boucle non perturb\'ee
(que l'on peut approcher par $\cal A$ puisque les boucles non
perturb\'ees longent les s\'eparatrices ), que le terme 
non perturb\'e de l'action le long d'une boucle perturb\'ee.
On se ram\`ene \`a montrer que 
${\cal F}(\gamma_{\omega_1}) \le {\cal F} (\partial B)$.
Soit $(T_0,Q_0)\in \partial B$. 
Posons 
\begin{eqnarray*}
      P = &  -& 
      \int _0 ^{\infty}
      (1-\cos q_0(t-T_{-1}))\,f(Q_{-1}+\omega_1(t-T_{-1}),q_0(t-T_{-1}),t)\,dt
        \\
       &  - & 
      \int _{-\infty} ^0
      (1-\cos q_0(t-T_1))\,f(Q_1+\omega_1(t-T_1),q_0(t-T_1),t)\,dt 
      .
\end{eqnarray*}
Comme indiqu\'e ci--dessus, 
les approximations suivantes se d\'emontrent
sans difficult\'e gr\^ace au lemme 1:
$${\cal P}(T_0,Q_0) = {\cal A} (\omega_1,T_0,Q_0) +P+ \epsilon(\mu),$$
$${\cal P}(\gamma_{\omega_1}) = 
  {\cal A} ({\omega_1},T_{\omega_1},Q_{\omega_1}) +P+ \epsilon(\mu),$$
 o\`u $\epsilon(\mu) \longrightarrow 0$  quand $\mu \longrightarrow 0$.
Donc, en utilisant la condition 1 :
$${\cal P}(T_0,Q_0)\ge 
  {\cal P}(\gamma_{\omega_1})+A+\epsilon(\mu).$$
Ceci d\'emontre la premi\`ere partie de nos assertions: le terme de 
perturbation $\cal P$ admet un minimum local sur $B$.
Posons: $\gamma_{\omega_1}(t)=(Q_m(t),q_m(t))$.
Les in\'egalit\'es suivantes expriment que, comme annonc\'e, les variations de
l'action non perturb\'ee sont plus faibles que celles de $\cal P$.
$$
  \int_{T_{-1}}^{T1} \left( \frac 1 2 \dot q_0 ^2 + (1-\cos q_0) \right) \ge
  \int_{T_{-1}}^{T1} \left( \frac 1 2 \dot q_m ^2 + (1-\cos q_m) \right)
  + o(\mu).
$$
On d\'emontre cette in\'egalit\'e en remarquant que, si $(Q(t),q(t))$ est 
la r\'eunion de boucles non perturb\'ees de m\^emes extr\'emit\'es que 
$(Q_0(t),q_0(t))$, alors $\int_{T_{-1}}^{T1} 1/2\,\dot q_0 ^2 + (1-\cos q_0) 
\ge
\int_{T_{-1}}^{T1} 1/2\,\dot q ^2 + (1-\cos q)$. On est alors ramen\'e \`a
comparer les actions de deux trajectoires du pendule. 
L'in\'egalit\'e
$$  \int_{T_{-1}}^{T1} \left( \frac 1 2 \dot Q_0 ^2 - 
    \frac 1 2 \dot Q_m ^2 \right) \,\ge -\frac C {A_0} \mu $$
permet de conclure:
\begin{eqnarray*}
  {\cal F}(T_0,Q_0)   
  & =
  & {\cal F}_0(T_0,Q_0)+\mu{\cal P}(T_0,Q_0)\\
  & \ge 
  & {\cal F}_0(\gamma_{\omega_1})+\mu{\cal P}(\gamma_{\omega_1})
  + o(\mu)+ A\mu - \frac C {A_0} \mu \\
  & \ge
  & {\cal F}(\gamma_{\omega_1})+\left( A-\frac C {A_0} \right)\mu
  +o(\mu).
\end{eqnarray*}
Si $A_0$ est suffisamment grand
(de sorte que $(A-C/A_0)>0$),  pour $\mu$
petit, $\cal F$ atteint son minimum sur l'int\'erieur de $B$. Le recollement
est donc effectu\'e : on a construit une transition simple, et donc 
d\'emontr\'e la proposition.

\end{document}